\documentclass[12pt,a4paper]{article}

\usepackage[left=1cm,right=1cm,top=2cm,bottom=2cm]{geometry}
\usepackage[dvips]{graphics}
\usepackage[dvips]{epsfig}
\usepackage{color}
\usepackage{hyperref}
\usepackage{tikz,authblk}

\newtheorem{lemma}{Lemma}
\newtheorem{theorem}[lemma]{Theorem}
\newtheorem{corollary}[lemma]{Corollary}
\newtheorem{proposition}[lemma]{Proposition}
\newtheorem{definition}{Definition}
\newtheorem{remark}{Remark}
\newtheorem{example}{Example}

\newcommand{\dimo}{\noindent \emph{Proof. }}
\newcommand{\qed}{\\ \rightline{$\Box$ \ \ \ \ \ \ \ \ \ \ \ \ \ \ \ }\\}
\newcommand{\e}{\varepsilon}
\newcommand{\G}{\Gamma}

\usepackage{latexsym}
\usepackage{amsfonts}
\usepackage{amssymb}
\usepackage{amsmath}
\usepackage[english]{babel}

\begin{document}

\title{Trisections of PL 4-manifolds arising from colored triangulations}

 \renewcommand{\Authfont}{\scshape\small}
 \renewcommand{\Affilfont}{\itshape\small}
 \renewcommand{\Authand}{ and }

\author[1] {Maria Rita Casali}
\author[2] {Paola Cristofori}

\affil[1,2] {Department of Physics, Informatics and Mathematics, University of Modena and Reggio Emilia\par  Via Campi 213 B, I-41125 Modena (Italy), casali@unimore.it,\ paola.cristofori@unimore.it}

\maketitle

\abstract{The purpose of the present paper is twofold: firstly to extend to non-orientable compact 4-manifolds the notion of {\it gem-induced trisection}, directly obtained from colored triangulations (or, equivalently, from colored graphs encoding them, called {\it gems}); secondly to prove that, both in the orientable and non-orientable case, if the boundary  is homeomorphic to a connected sum of sphere bundles over $\mathbb S^1$,  gem-induced trisections naturally give rise to trisections of the corresponding closed 4-manifold. 
As a consequence, an estimation of the trisection genus of any closed orientable $4$-manifold in terms of surgery description  is obtained via  colored triangulations. }

\endabstract
 
\bigskip
  \par \noindent
  {\small {\bf Keywords}: trisection, compact 4-manifold, colored triangulation}

 \medskip
  \noindent {\small {\bf 2020 Mathematics Subject Classification}: 57Q15 - 57K40 - 57M15}

\maketitle

\section{Introduction}
\label{intro}

In 2018 (see \cite{Bell-et-al}), Bell, Hass, Rubinstein and Tillmann  introduced the use of (singular) triangulations to study {\it trisections} of closed orientable $4$-manifolds, i.e. decompositions into three $4$-dimensional handlebodies, mutually intersecting in $3$-dimensional handlebodies and globally intersecting in a closed surface  (according to the definition by Gay and Kirby in \cite{Gay-Kirby}).   
The same combinatorial approach has been applied in \cite{SpreerTillmann(Exp)} to the special case of {\it colored triangulations} 
of (closed) PL $4$-manifolds with exactly ten edges, in order to compute the {\it trisection genus} (i.e. the minimum genus of the intersecting surface)  of connected sums of $\mathbb{CP}^2$ - possibly with reversed orientation -, $\mathbb S^2\times\mathbb S^2$ and the $K3$-surface.  

The above ideas have been generalized in  \cite{Casali-Cristofori gem-induced}, where the assumption on the number of edges of the colored triangulation is relaxed and orientable 4-manifolds with connected boundary are included, too: the resulting decompositions are called {\it gem-induced trisections}, since they are directly obtained from {\it gems} (i.e. {\it Graphs Encoding Manifolds}). 

While in the closed setting a gem-induced trisection is a particular kind of trisection, when the boundary is non-empty  the decomposition is different  from that introduced in  \cite{Castro-Gay-Pinzon}, which is commonly intended as a trisection in the boundary case: in fact, it doesn't involve compression bodies and open book decompositions, but only handlebodies and Heegaard splittings, according to a suggestion by Rubinstein and Tillmann in  \cite[Section 4.5]{Rubinstein-Tillmann}.

On the other hand, in  \cite{Miller-Naylor}, the classical definition of trisection has been broadened to non-orientable manifolds.   
It is then natural to address the problem of  extending to the non-orientable case the notion of gem-induced trisection of a compact 4-manifold.   

Section \ref{s.gem-induced trisections} of the present paper describes how such a decomposition can be obtained from a colored triangulation of the associated singular manifold, or, equivalently, from its dual graph, which is called a {\it gem} of the compact manifold itself (see Section \ref{preliminaries} for a brief overview of this combinatorial representation for compact PL manifolds).  

\medskip 

Actually, gem-induced trisections turn out to fit the definition of trisection  - hence allowing a direct estimation of the trisection genus - only for a suitable class of  closed orientable $4$-manifolds (possibly comprehending all simply-connected ones, according to Kirby problem n. 50):   see Proposition \ref{trisection_vs_gem-induced}.

Nevertheless, the particular type of extension to the boundary case performed by gem-induced trisections allows also an ``indirect'' approach to trisections of closed (orientable  and  non-orientable) 4-manifolds.    In fact, in Section \ref{s.gem-induced vs trisections} we prove that, in case of a compact 4-manifold  $M$ whose boundary is homeomorphic to a connected sum of sphere bundles over $\mathbb S^1$, any  gem-induced trisection of $M$ naturally gives rise to a trisection of the associated closed 4-manifold $\bar M$ with the same intersecting surface; 
therefore, an upper bound for the trisection genus of $\bar M$ is obtained (see Theorem \ref{trisection_from_gem-induced}). 

\medskip 

 In particular, trisections arising from colored triangulations turn out to minimize the trisection genus for a wide class of (orientable and non-orientable) $4$-manifolds: see Proposition  \ref{calculations}.   
 
 \bigskip

Moreover,  combining the above results with those contained in  \cite{Casali-Cristofori Kirby-diagrams} and \cite{Casali-Cristofori gem-induced},  we prove the existence, for each closed orientable $4$-manifold $\bar M$, of trisections  arising from (suitable) gems of the compact 4-manifold $M$, that consists of the 0-, 1- and 2-handles in a handle-decomposition of $\bar M$.
Hence, an estimation of the trisection genus $g_T(\bar M)$  in terms of surgery description is obtained (see Section \ref{s.gem-induced vs trisections} for details): 

\begin{theorem}\label{trisection_from_Kirby-diagram} \ 
\begin{itemize}
\item[(i)]
 For each closed orientable $4$-manifold $\bar M$,  
$$g_T(\bar M) \le s+1,$$  $s$ being the crossing number  of a  (connected and with dotted components in ``good position'') Kirby diagram representing $\bar M$. 
\item[(ii)] Furthermore, if  $\bar M$ admits a handle decomposition lacking in $1$-handles, then 
$$g_T(\bar M) \le m_\alpha, $$
$m_{\alpha}$ being  the number of $\alpha$-colored regions in a chess-board coloration of a  (connected and with no dotted component) Kirby diagram representing  $\bar M$.
\end{itemize}
\end{theorem}

\section{Colored triangulations and gems of compact PL ma\-ni\-folds} 
\label{preliminaries}
Throughout this paper we will work in the PL category, therefore manifolds and maps under consideration will always be PL.

By a {\em singular $n$-manifold} we mean a closed connected $n$-dimensional polyhedron admitting a simplicial triangulation where the link of vertices are closed connected $(n-1)$-manifolds while the links of the $h$-simplices, with $h>0$, are $(n-h-1)$-spheres. Vertices whose links are not spheres will be called {\em singular}.

A {\em colored triangulation} of a singular $n$-manifold $N$ is a (pseudo)trian\-gu\-lation K of $N$ endowed with a labeling of its vertices by the set of integers $\Delta_n=\{0,\ldots,n\}$ which is injective on each simplex.

\smallskip

Such a colored triangulation is usually combinatorially visualized by the 1-skeleton $\G(K)$ of its dual complex, endowed with the edge-coloration inherited from the labeling of $K$: an edge $e$ of $\G(K)$ has color $c\in\Delta_n$ iff  no vertex of the $(n-1)$-simplex of $K$ dual to $e$ is labeled $c.$  
$\G(K)$ is an {\it $(n+1)$-colored graph} (i.e. a regular multigraph with degree $n+1$, such that its edge-coloration by $\Delta_n$ is injective on adjacent edges\footnote{According to a well-established literature (see, for example,  \cite{Casali-Cristofori-Gagliardi RIMUT 2020} and references within), an $n$-dimensional pseudocomplex $K(\G)$ may be associated to any $(n+1)$-colored graph $\G$, by considering an $n$-simplex, with vertices labelled by the elements of $\Delta_n,$ for each vertex of $\G$ and  by gluing two $n$-simplices  along their $(n-1)$-dimensional faces opposite to $c$-labeled vertices, whenever the corresponding vertices of $\G$ are $c$-adjacent ($c\in\Delta_n$). Note that, if $|K(\G)|$ is a singular manifold, both $\G(K(\G))=\G$ and $K(\G(K))= K$ hold, so colored triangulations of singular manifolds and colored graphs representing them can be associated unambiguously.}), which is said to {\em represent} the singular manifold $N.$

\begin{remark} {\em 
The duality between $K$ and $\G(K)$ establishes a bijective correspondence  between the 
$(n-h)$-simplices of $K$ whose vertices are labeled by  $\Delta_n-\{c_1, \dots, c_h\}$ 
and the  connected components of the subgraph $\G_{\{c_1, \dots, c_h\}}$ obtained from $\G(K)$  by 
considering only edges colored by $\{c_1, \dots, c_h\}$, which are called {\it $\{c_1, \dots, c_h\}$-residues}. 
In particular, the connected components obtained from $\G(K)$ by deleting all $c$-colored edges ($c\in\Delta_n$), called {\it $\hat c$-residues} of $\G(K)$, are $n$-colored graphs representing the disjoint links\footnote{Given an $h$-simplex $\sigma^h$ of $K$, the {\it disjoint star} of $\sigma^h$ in $K$ is the pseudocomplex obtained by taking all $n$-simplices of $K$ having $\sigma^h$ as a face and re-identifying only their faces that contain $\sigma^h.$ The {\it disjoint link} of $\sigma^h$ in $K$ is the subcomplex of the disjoint star formed by those simplices that do not intersect $\sigma^h.$} of the $c$-labeled vertices of $K$ (which are closed $(n-1)$-manifolds, since $K$ triangulates a singular $n$-manifold).   

\noindent As a consequence, by setting from now on $\G=\G(K)$, the following characterization holds:

\centerline{\it $\vert K\vert$ is a closed $n$-manifold iff, for each color $c\in\Delta_n$,}
\centerline{\it all $\hat c$-residues of $\G$ represent the $(n-1)$-sphere.}  

\smallskip

A {\it $\hat c$-residue} of $\G$ will be called {\it singular} if it corresponds to a singular vertex of $\vert K \vert$, i.e. if it does not represent the sphere. A color $c$ is said to be {\it singular} if at least one $\hat c$-residue is singular. 

Moreover, the number of $\{c_1, \dots, c_h\}$-residues will be denoted by $g_{c_1, \dots, c_h}$ (or, for short, by $g_{\hat c}$ if $h = n$ and $\{c \}= \Delta_n - \{c_1, \dots, c_n\}$). }
\end{remark}  

\medskip   

If $M$ is a  compact $n$-manifold, then a singular $n$-manifold $\widehat M$ can be constructed by capping off each component of $\partial M$ by a cone over it. 
Therefore, $(n+1)$--colored graphs may be used to represent compact PL $n$-manifolds, as well:  
 
\begin{definition}\label{def.gem} {\em
An $(n+1)$-colored graph {\it represents} a compact $n$-manifold $M$ (or, equivalently, it is a {\it gem} of $M$, where gem means {\it Graph Encoding Manifold})  if and only if  the dual pseudocomplex is a colored triangulation of the singular manifold $\widehat M$.}
\end{definition} 

The so called {\it gem theory}, or {\it crystallization theory}, is based on the following existence result (which extends the one originally stated in \cite{Pezzana}  for closed manifolds):  

\begin{theorem}{\em \cite{Casali-Cristofori-Grasselli}}\ \label{Theorem_gem}  
Any compact orientable (resp. non orientable) $n$-manifold $M$
admits a bipartite (resp. non-bipartite) $(n+1)$-colored graph $\G$ representing it.
\par \noindent In particular, if  $M$ has empty or connected boundary: 
\begin{itemize}
\item $\G$ may be assumed to have color $n$ as its unique possible singular color, and exactly one $\hat n$-residue 
(in this case we will say that $\G$ belongs to the class $G^{(n)}_s$);\footnote{Equivalently, the associated colored triangulation $K$ of $\widehat M$  has exactly one vertex labelled $n$,  which is the only possible singular vertex of $K$ 
(in this case we will say that $K$ belongs to the class $\mathcal K^{(n)}_s$).
Note that, obviously, the vertex of $K$ labelled $n$ (or, equivalently, the $\hat n$-residue of $\G$) is singular if and only if $\partial M$ is non-empty and non-spherical.} 
\item $\G$ may be assumed to have exactly one $\hat c$-residue, $\forall c \in \Delta_n$ (in this case $\G$ is  called  
a \emph{crystallization} of $M$).\footnote{Equivalently, the associated colored triangulation of $\widehat M$   has exactly $n+1$ vertices.}  
\end{itemize}
\end{theorem}

\smallskip

In this paper only manifolds with empty or connected boundary will be considered. Note that, according to Definition \ref{def.gem}, the same $(n+1)$-colored graph can represent both a closed $n$-manifold and the compact $n$-manifold  obtained by deleting from it the interior of an $n$-ball. In order to avoid this ambiguity, from now on we will also exclude manifolds with spherical boundary. 
However, all definitions and results of this paper regarding closed manifolds (for which, obviously, $M=\widehat M$ holds) can be easily translated to fit the case of spherical boundary.

\smallskip

Let us recall (\cite{Gagliardi 1981}) that, given a connected bipartite (resp. non-bipartite)  
$(n+1)$-colored graph $\G$ of order $2p$, then for each cyclic permutation $\varepsilon = (\varepsilon_0,\ldots,\varepsilon_n)$ of $\Delta_n$, up to inverse, there exists a particular type of cellular embedding, called  \emph{regular}\footnote{More precisely, a regular embedding is a cellular embedding whose regions are bounded by the images of the $\{\varepsilon_j,\varepsilon_{j+1}\}$-colored cycles of $\G$, for each $j \in \mathbb Z_{n+1}$.}, of $\G$  into an orientable (resp. non-orientable) closed surface. 
Moreover, the genus (resp. half the genus) of this surface, denoted by $\rho_{\varepsilon} (\G)$, satisfies
\begin{equation} \label{eq.regulargenus} 
2 - 2\rho_\varepsilon(\G)= \sum_{j\in \mathbb{Z}_{n+1}} g_{\varepsilon_j, \varepsilon_{j+1}} + (1-n)p.
\end{equation}
\smallskip

The {\em regular genus} of  $\G$ is defined as
$$\rho(\G) = min\{\rho_\varepsilon(\G)\ \vert \ \varepsilon\ \text{is a cyclic permutation of \ } \Delta_n\}. $$

The  {\it regular genus} of a compact $n$-manifold $M$ is a PL invariant extending to higher dimension the classical genus of a surface and the Heegaard genus of a $3$-manifold: it is defined as
$$\mathcal G (M) = min\{\rho(\G)\ \vert \ \G\ \text{is a gem of \ } M\}.$$  

It was proved  in \cite{Ferri-Gagliardi Proc AMS 1982} that regular genus zero characterizes spheres in any dimension; morever,  other classification results via regular genus are available, especially in dimension $4$ and $5$ (see \cite{Casali-Gagliardi ProcAMS},  \cite{Casali_Forum2003},  \cite{generalized-genus}  and their references). Sections \ref{s.gem-induced trisections} and  \ref{s.gem-induced vs trisections} will show that, in dimension $4$, the regular genus is also strongly involved in the estimation of trisection genus via colored triangulations.

\section{Gem-induced trisections of non-orientable 4-manifolds}  \label{s.gem-induced trisections}

As already recalled in Section 1, gem-induced trisections of compact orientable 4-manifolds with empty or connected boundary have been introduced in \cite{Casali-Cristofori gem-induced}, as an effective tool for the study of trisections in the PL setting. 
They have been also used in \cite{Martini-Toriumi} to investigate connections between trisections of closed orientable 4-manifolds and {\it colored tensor models}  (\cite{Casali-Cristofori-Dartois-Grasselli}). 

In this section we will show how the concept of gem-induced trisection can be extended so as to include also the non-orientable case.
\smallskip

Let $M$ be a compact 4-manifold with empty or connected boundary and let $K$ be a colored triangulation of the singular manifold $\widehat M$. Note that, by Theorem \ref{Theorem_gem}, it is always possible to suppose $K\in\mathcal K^{(4)}_s,$ i.e. to assume that $K$ has only one possible singular vertex, labelled by color 4.
Furthermore, let us fix a cyclic permutation $\e$ of $\Delta_4$ and suppose, for simplicity, $\e_4=4.$

\begin{figure}[t]
\centering
\scalebox{0.35}{\includegraphics{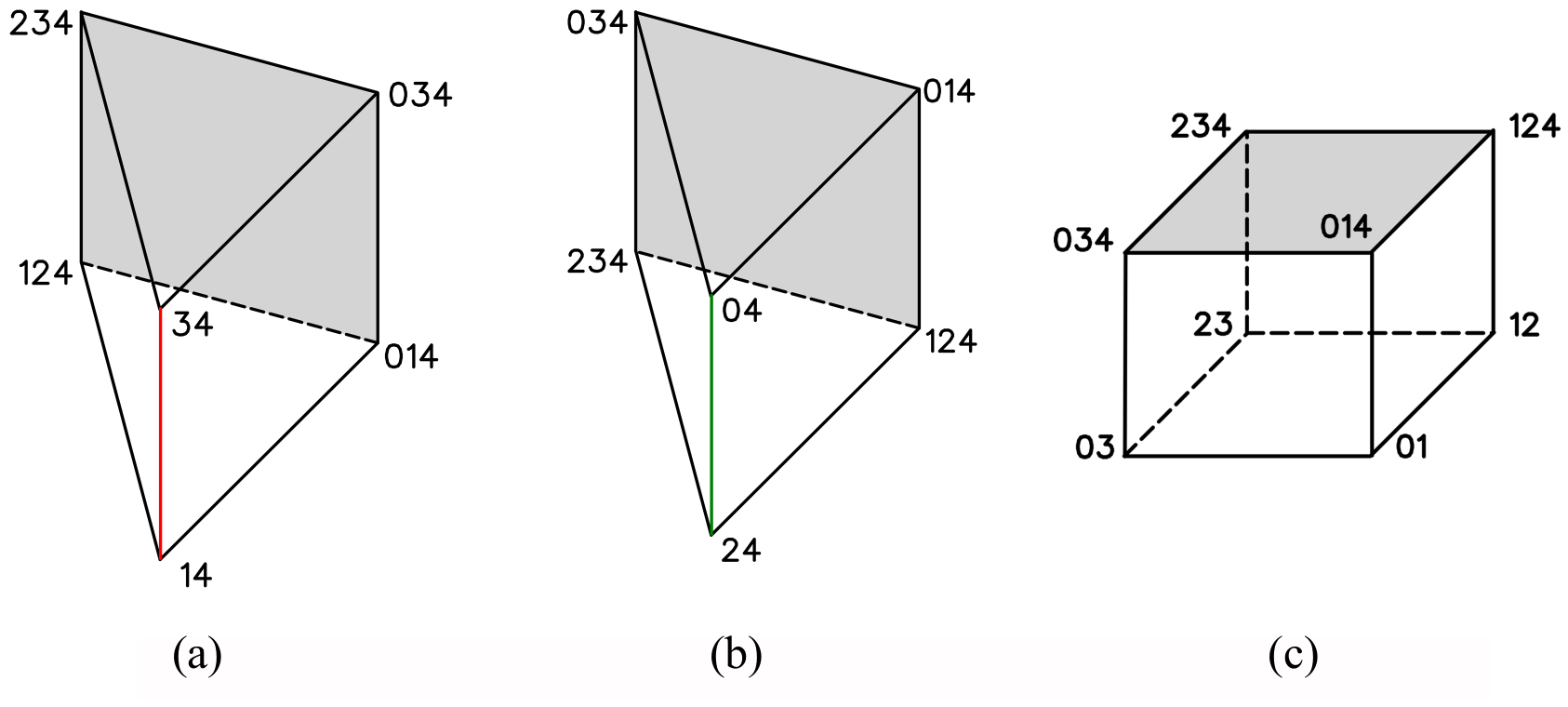}}
\caption{the intersections of a 4-simplex of $K(\G)$ with $H_{01}$,  $H_{02}$ and $H_{12}$ respectively}   
\label{fig.prism-cube}
\end{figure}

Then - by generalizing an idea of \cite{Bell-et-al} and \cite{SpreerTillmann(Exp)} -  we can construct a particular decomposition of the singular manifold $\widehat M$ associated to $K$  
and $\e$ in the following way:

\begin{itemize}
\item let $\sigma$ be the standard $2$-simplex, whose vertices are denoted by $v_0, v_1, v_2$ and let $\mu\ :\ K(\G)\to\sigma$ be the simplicial map defined by $\mu(v)=v_1$ (resp. $\mu(v)=v_2$) iff $v$ is a vertex of $K$ 
labelled by $\e_0$ or $\e_2$ (resp. $\e_1$ or $\e_3$); finally let us set $\mu(w)=v_0$, where $w$ is the unique $\e_4$-labelled vertex of $K.$   

\item The preimage $\widehat H_0$ of the star of $v_0$ in the first barycentric subdivision $\sigma^\prime$ of $\sigma$ is the cone over the disjoint link of $w$, while the preimage  $H_1$ (resp. $H_2$) of the star of $v_{1}$ (resp. $v_{2}$) in  $\sigma^\prime$ is a regular neighbourhood of the 1-subcomplex of $K$  generated by the $\e_0$- and $\e_2$-labelled (resp. $\e_1$- and $\e_3$-labelled) vertices. Therefore both $H_1$ and $H_2$ are 4-dimensional handlebodies.

\item The 3-dimensional subcomplex $H_{01}=\widehat H_0\cap H_1$ (resp. $H_{02}=\widehat H_0\cap H_2$) is a 3-dimensional handlebody. In fact,  its intersection with each 4-simplex of $K$  is a triangular prism as the one in Figure \ref{fig.prism-cube}(a)  (resp. Figure \ref{fig.prism-cube}(b)), whose vertices are barycenters of a 1- or 2-simplex,  which we indicate by the labels of its spanning vertices. The quadrangular face of the prism opposite to the blue (resp. red) edge is free; so, $H_{01}$ (resp. $H_{02}$) collapses to a graph.

\item Note that  $\Sigma=\widehat H_0\cap H_1\cap H_2$  is  formed by the free faces of the above prisms (see  Figure \ref{fig.prism-cube}, where $H_{12}=H_1\cap H_2$), one for each 4-simplex of $K$;  since each edge of such quadrangular faces corresponds to a 3-simplex of $K$, it is shared by exactly two faces and hence $\Sigma$ turns out to be  a closed connected surface. Moreover,  its Euler characteristic is $\chi(\Sigma)=g_{0,1}+g_{1,2}+g_{2,3}+g_{0,3}-4p+2p =  2 - 2\rho_{\varepsilon_{\hat 4}}(\G(K)_{\hat 4})$, where $ \varepsilon_{\hat 4}=(\e_0,\e_1,\e_2,\e_3).$
\end{itemize}

As a consequence of the above construction, the following result can be stated, which extends to the non-orientable case the analogous one in \cite{Casali-Cristofori gem-induced}.

\begin{proposition}\label{genus_trisection} Let $\widehat M$ be a singular $4$-manifold with one singular vertex at most. For each colored triangulation $K \in \mathcal K^{(4)}_s$ of $\widehat M$ and for each cyclic permutation $\e=(\e_0,\e_1,\e_2,\e_3,4)$ of $\Delta_4$, the triple $(\widehat H_{0}, H_{1},H_{2})$ satisfies the following properties: 
\begin{itemize}
 \item [(i)] $\widehat M = \widehat H_{0}\cup H_{1}\cup H_{2}$ and the interiors of $\widehat H_{0},$ $H_{1},$ $H_{2}$ are pairwise disjoint;
\item [(ii)] both $H_1$ and $H_2$ are 4-dimensional handlebodies of genus $g_{\e_1,\e_3,\e_4} - g_{\widehat{\e_0}} - g_{\widehat{\e_2}} + 1$  and $g_{\e_0,\e_2,\e_4} -  g_{\widehat{\e_1}} - g_{\widehat{\e_3}} + 1$ respectively;
 \item [(iii)] $\widehat H_{0}$ is homeomorphic to the cone over the disjoint link of the singular vertex of $\widehat M$ (or to a $4$-ball, if   $\widehat M$ is a $4$-manifold); 
 \item [(iv)] $H_{01}=\widehat H_{0}\cap H_{1}$ and $H_{02}=\widehat H_{0}\cap H_{2}$ are $3$-dimensional handlebodies; 
\item  [(v)]  $\Sigma=\widehat H_{0}\cap H_{1}\cap H_{2}$ is a closed connected surface. 
\end{itemize}
Moreover, if $H_{12}=H_1\cap H_2$ is a 3-dimensional handlebody, too, then all the above handlebodies, as well as the surface $\Sigma$, are orientable or not according to the orientability of $\widehat M.$ Therefore, in the first case $\Sigma$ has genus $\rho_{\varepsilon_{\hat 4}}(\G(K)_{\hat 4})$, while in the second one it has genus $2\rho_{\varepsilon_{\hat 4}}(\G(K)_{\hat 4}).$
\end{proposition}

\dimo       
 Statements (i), (iii), (iv) and (v) directly follow by construction, as well as the fact that $H_{1}$ and $H_{2}$ are handlebodies. To complete the proof of statement (ii), note that the $1$-dimensional subcomplex $K_{\e_0\e_2}$ (resp. $K_{\e_1\e_3}$) of $K$ generated by the $\{\e_0,\e_2\}$-colored (resp. $\{\e_1,\e_3\}$-colored) vertices, has exactly $g_{\e_1,\e_3,4}$  (resp. $g_{\e_0,\e_2,4}$) edges and $g_{\widehat{\e_0}} + g_{\widehat{\e_2}}$ (resp. $g_{\widehat{\e_1}} + g_{\widehat{\e_3}}$) vertices.
Since $K$ is a pseudomanifold, $K_{\e_0\e_2}$ (resp. $K_{\e_1\e_3}$) is connected; hence $H_1$ (resp. $H_2$) has genus  $g_{\e_1,\e_3,4} - g_{\widehat{\e_0}} - g_{\widehat{\e_2}} + 1$ (resp. $g_{\e_0,\e_2,4} - g_{\widehat{\e_1}} - g_{\widehat{\e_3}} + 1$).

\smallskip  
 
Let us now observe that $(H_{01},H_{02},\Sigma)$ is a Heegaard splitting of $\partial\widehat H_0$, which is  the disjoint link of the (possibly singular) vertex $w$ of $K$; moreover, if  $H_{01},H_{02},H_{12}$ are all handlebodies, then for each $i\in \{1,2\}$,  $(H_{ij},H_{ik},\Sigma)$, with $\{j,k\}= \{0,1,2\}-\{i\}$, is  a Heegaard splitting of the 3-manifold $\partial H_i$,  which is a connected sum of  copies of the $\mathbb S^2$-bundle over $\mathbb S^1$,  orientable or  non-orientable according to $H_i$.
Therefore, it is obvious that, if $\widehat M$ is orientable,  all 4-dimensional ``pieces'' $\widehat H_0,$ $H_1$ and $H_2$  are orientable and the same holds for the 3-dimensional handlebodies and the surface $\Sigma.$

On the other hand, if $\widehat M$ is non-orientable, then all the handlebodies and the surface $\Sigma$ must be non-orientable. In fact, the existence of the above Heegaard splittings allows easily to check that, if one of the 4-dimensional ``pieces'' were orientable, then the Heegaard splitting of its boundary would be formed by orientable elements  and, as a consequence, both the third  3-dimensional handlebody and the other two 4-dimensional ``pieces'' would be orientable, too.

In the same way, if one of the 3-dimensional handlebodies, say $H_{ij}$, (resp. if the surface $\Sigma$) were orientable, the existence of the above Heegaard splittings would imply the orientability of the other 3-dimensional handlebodies, as well as of the boundaries of both  the 4-dimensional ``pieces'' intersecting in $H_{ij}$ (resp. of all 3-dimensional handlebodies, as well as of the boundaries of all 4-dimensional ``pieces''), and therefore of all 4-dimensional ``pieces'', too. 

Moreover, the final statement regarding the genus of $\Sigma$ is a trivial consequence of the computation  $\chi(\Sigma)= 2 - 2\rho_{\varepsilon_{\hat 4}}(\G(K)_{\hat 4}).$    
\qed       

\smallskip

It is easy to check that, if $\widehat M$ is the singular manifold associated to a compact 4-manifold $M$, the above triple $(\widehat H_{0},H_{1},H_{2})$ naturally induces a decomposition $\mathcal T(\Gamma(K),\e)=(H_0,H_1,H_2)$ of $M$, where  $H_0$ is the collar on $\partial M$ obtained by deleting from $\widehat H_0$ a suitable neighbourhood of the singular vertex of $K$ (or $H_0=\widehat H_0$ is a $4$-ball, if  $\partial M=\emptyset$).  
   
\smallskip

In full analogy with what was already introduced in the orientable case, we can give the following definition:

\begin{definition}\label{def.gem-induced} {\em 
Let $M$ be a compact 4-manifold with empty or connected boundary.  A {\em gem-induced trisection} of $M$ is a decomposition $\mathcal T(\G(K),\e)=(H_0,H_1,H_2)$ of $M$ such that $H_{12}$ is a 3-dimensional handlebody, $K \in \mathcal K^{(4)}_s$ being a colored triangulation of $\widehat M$ and $\e=(\e_0,\e_1,\e_2,\e_3,4)$ a cyclic permutation of $\Delta_4$.
\par \noindent 
In this case, $\Sigma$ is called the {\em central surface} of the gem-induced trisection, and we refer to the common genus of all 3-dimensional handlebodies as the {\em genus} of the gem-induced trisection.   }
\end{definition}

The term ``gem-induced'' refers to the fact that the obtained decomposition can be thought as induced both by the triangulation $K$ of $\widehat M$ and by its dual 5-colored graph  $\G(K)\in G_s^{(4)}$, which in turn is also a gem of the manifold $M.$   
Moreover, in the following we will often refer explicitly only to the graph, by using the notations $\mathcal T(\Gamma,\e)$ for the gem-induced trisection and $genus(\mathcal T(\G,\e))$ for its genus (which is equal to $\rho_{\varepsilon_{\hat 4}}(\G_{\hat 4})$).     

\begin{remark}\label{remark-Heegaard} 
{\em 
Note that, in the case of a gem-induced trisection $\mathcal T(\G,\e)=(H_0,H_1,H_2)$ of a compact 4-manifold $M$ with empty or connected boundary, then $(H_{01},H_{02},\Sigma)$ is a Heegaard splitting of $\partial\widehat H_0 = \partial H_0 =\ \mathbb S^3$  if $M$ is closed,  while, if $\partial M\neq\emptyset$, it is a Heegaard splitting of $\partial\widehat H_0 = \partial M.$ 

Moreover, if $M$ (and so $\widehat M$) is non-orientable, $\partial\widehat H_0$ must be non-orientable; as a consequence no closed non-orientable 4-manifold can admit a gem-induced trisection and the same happens when $M$ is non-orientable with orientable boundary.  }
\end{remark}

\begin{remark}\label{remark-simplyconnected} {\em 
Let $\mathcal T(\G,\e)=(H_0,H_1,H_2)$ be a gem-induced trisection of a compact 
4-manifold $M$ with empty or connected boundary. 
The existence of the Heegaard splittings of $\partial H_1$ and $\partial H_2$ given by the 3-dimensional handlebodies of the decomposition, together with  suitable applications of Van-Kampen's theorem, implies that $\pi_1(\widehat H_0)$ surjects to $\pi_1(\widehat M)$ (details can be found in \cite{Casali-Cristofori gem-induced}, since the arguments do not depend on orientability).  Since $\widehat H_0$ is contractible, the simply-connectedness of $\widehat M$ is a necessary condition for the existence of gem-induced trisections of $M.$ }
\end{remark}

By construction and by Definition \ref{def.gem-induced}, a gem-induced trisection of a closed $4$-manifold  is a particular type of trisection (in the sense of  \cite{Gay-Kirby}; see also Section \ref{s.gem-induced vs trisections}), where one of the $4$-dimensional handlebodies is a $4$-ball. 
On the other hand, in the boundary case, gem-induced trisections provide a decomposition which differs in substance from the extended notion of trisection introduced in \cite{Castro-Gay-Pinzon} (see \cite[Remark 16]{Casali-Cristofori gem-induced}).
Therefore, as a consequence of Remark \ref{remark-Heegaard} and \cite[Proposition 25]{Casali-Cristofori gem-induced}, which characterizes closed orientable $4$-manifolds admitting gem-induced trisections, we  can state that gem-induced trisections fit the usual definition of trisection only for a particular class of closed orientable $4$-manifolds; however, it is worthwhile to note that this class possibly comprehends all simply-connected ones, according to Kirby problem n. 50.

\begin{proposition}\label{trisection_vs_gem-induced} 
Let $M$ be a compact PL $4$-manifold with empty or connected boundary. 
A gem-induced trisection of $M$ is a trisection if and only if  $M$ is closed and orientable. 

\noindent Moreover, a closed orientable $4$-manifold admits a gem-induced trisection if and only if it admits a handle decomposition lacking in 3-handles (or in 1-handles).
\end{proposition}
\ \qed

\medskip

In the orientable setting, a combinatorial condition ensuring  $\mathcal T(\G,\e)$ to be a gem-induced trisection of $M$ was already presented in \cite{Casali-Cristofori gem-induced}; the same condition works as well also in the non-orientable case: 
  
\begin{proposition} \label{CS gem-induced trisections}
Let $M$ be a compact $4$-manifold with empty or connected boundary and $\G\in G_s^{(4)}$ a gem of $M$ of order $2p$;  if  there exists an ordering $(e_1,\ldots,e_p)$ of the 4-colored edges of $\G$ such that for each $j\in\{1,\ldots,p\}$:
\centerline{(*)\hskip 25pt  there exists $i\in\Delta_3$  such that all 4-colored edges of  the \hskip 35pt} 
\centerline{$\{4,i\}$-colored cycle containing $e_j$ belong to the set $\{e_1,\ldots,e_j\}$,}
then $\mathcal T(\G, \varepsilon)$ is a gem-induced trisection of $M$, for each cyclic permutation $\e$ of $\Delta_4.$       
\end{proposition} 
 
\dimo      
The involved arguments are exactly the same as in the orientable case: see the proof of \cite[Proposition 20]{Casali-Cristofori gem-induced}, together with the final part of the proof of \cite[Proposition 17]{Casali-Cristofori gem-induced}.         
\qed

\begin{definition}\label{def_GT-genus} {\em 
Let $M$ be a compact $4$-manifold with empty or connected  boundary, that admits gem-induced trisections. The  \emph{G-trisection genus} % $g_{GT}(M)$ 
of $M$ is defined as: 
$$   g_{GT}(M) \, = \,  \min \{genus(\mathcal T (\G, \e))  \  |  \ \mathcal T (\G, \e)  \ \ \text{is a gem-induced trisection of }  M\}.$$  }
\end{definition}

\medskip

\medskip

The properties of the G-trisection genus which were already proved in  \cite{Casali-Cristofori gem-induced} for orientable 4-manifolds can be easily extended to the non-orientable case, since the arguments do not depend on orientability.

\begin{proposition}\label{GTgenus-properties} Let $M$ be a compact 4-manifold with empty or connected boundary. Then:       
\begin{itemize}  
\item[(i)] $g_{GT}(M)\leq \rho_{\varepsilon_{\hat 4}}(\G_{\hat 4}),$ \  for any gem $\G$ of $M$ so that $ \mathcal T (\G, \e)$ is a gem-induced trisection.    
\item[(ii)]  $g_{GT}(M)=0 \ \ \Leftrightarrow  \ \ M \cong \mathbb S^4.$ 
\item[(iii)] If $M$ has non-empty boundary,  then \ $g_{GT}(M) \ge \mathcal H(\partial M),$ 
where $\mathcal  H(\partial M)$ denotes the Heegaard genus of the boundary.
\item[(iv)] The G-trisection genus is subadditive with respect to both the (internal) connected sum and the boundary connected sum.   
 \end{itemize} 
\end{proposition} 
\vskip-0.4truecm 
\ \qed 

\begin{figure}[t]
\centering
\scalebox{0.7}{\includegraphics{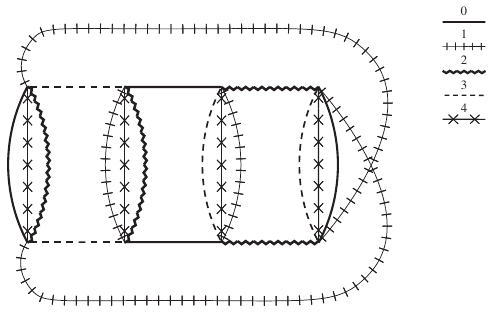}}
\caption{a gem of the genus one non-orientable handlebody $\tilde{\mathbb Y}^4_1$} 
\label{Y4_1-no}
\end{figure}

\begin{example}\label{example-handlebodies}
{\em 
Figure \ref{Y4_1-no} shows a gem $\G\in G_s^{(4)}$ of the genus one non-orientable handlebody $\tilde{\mathbb Y}^4_1.$ It is easy to see that $\mathcal T(\G,\e)$ is a gem-induced trisection of $\tilde{\mathbb Y}^4_1$ for each cyclic permutation $\e$ of $\Delta_4$, since $\G$ satisfies condition (*) 
of Proposition \ref{CS gem-induced trisections}; moreover, $\rho_{\varepsilon_{\hat 4}}(\G_{\hat 4})=1.$ Therefore, by statements (i) and (ii) of Proposition \ref{GTgenus-properties}, \ $ g_{GT}(\tilde{\mathbb Y}^4_1)=1.$

\noindent Furthermore, if $\tilde{\mathbb Y}^4_m$ denotes the genus $m$ non-orientable handlebody ($m \ge 1$), $g_{GT}(\tilde{\mathbb Y}^4_m)=m$ directly follows by statements  (iii) and (iv)  of Proposition \ref{GTgenus-properties}. The same value was already known from \cite{Casali-Cristofori gem-induced} to be the G-trisection genus of the genus $m$ orientable handlebody ${\mathbb Y}^4_m$. }
\end{example}

\section{Trisections arising from gem-induced trisections}  
\label{s.gem-induced vs trisections}

In the present section, we will show that gem-induced trisections of 4-manifolds with boundary allow an ``indirect'' approach to trisections of closed 4-manifolds, both in the orientable and non-orientable case. 

Let us first recall that, according to \cite{Gay-Kirby} and \cite{Miller-Naylor}, a  {\it trisection} of or a closed $4$-manifold $\bar M$  is a decomposition of $\bar M$ into three 4-dimensional handlebodies with disjoint interiors, whose pairwise intersections are 3-dimensional handlebodies  and all intersecting into a closed connected surface (called the {\it central surface} of the trisection).
Moreover, the seven ``pieces'' of the trisection (three in dimension 4, three in dimension 3, plus the central surface) are all orientable if and only if $\bar M$ is orientable and all non-orientable  if and only if $\bar M$ is non-orientable. It is also a direct consequence of the definition that all the 3-dimensional handlebodies have the same genus, which is called the {\it genus} of the trisection. 

The {\it trisection genus} $g_T(\bar M)$ of  any closed $4$-manifold $\bar M$ is defined as the minimum genus among all trisections of $\bar M.$ 

Throughout this section, the following lower bound for $g_T(\bar M)$ will turn out to be useful:  
\begin{equation}\label{genus_inequality} g_T(\bar M)\geq \beta_1(\bar M;\mathbb Z_2) + \beta_2(\bar M;\mathbb Z_2)
\end{equation}
where $\beta_i(\bar M;\mathbb Z_2)$ denotes the $i$-th Betti number of $\bar M$ with coefficient in $\mathbb Z_2.$

\smallskip

The proof of inequality \eqref{genus_inequality} is due to \cite{SpreerTillmann(Exp)},  where - however -  the trisection genus is defined as the minimal genus of a central surface, i.e. it coincides with  $g_T(\bar M)$ (resp.  $2g_T(\bar M)$) if $\bar M$ is orientable (resp. non-orientable). 

\medskip

The following theorem, by proving how trisections of closed $4$-manifolds may arise from gem-induced trisections of bounded $4$-manifolds, also yields an upper bound for the trisection genus in terms of the combinatorial invariant G-trisection genus defined in Section  \ref{s.gem-induced trisections}.    

 \begin{theorem}\label{trisection_from_gem-induced} 
\ \ Let \ $M$ \ be \ a \ compact \ \ orientable \ \ (resp. \ non-orientable) \ \ $4$-manifold \ with \ boundary  \  \  $\partial M \cong \#_m(\mathbb S^2 \times \mathbb S^1)$  (resp. \ $\partial M \cong \#_m (\mathbb S^2 \tilde \times \mathbb S^1$)), $m>0$.           
If $M$ admits a gem-induced trisection,  then the closed 4-manifold $\bar M$, uniquely obtained by gluing a $4$-dimensional handlebody along $\partial M$ (i.e. $\bar M \cong M \cup \mathbb Y^4_m$ or $\bar M \cong M \cup  \tilde {\mathbb Y}^4_m$, according to $M$ being orientable or not) admits a trisection with  the same central surface.  

As a consequence, $$g_T(\bar M) \le g_{GT} (M).$$ 
\end{theorem}
 
\dimo       
By hypothesis, a colored triangulation $K \in \mathcal K^{(4)}_s$ of  $\widehat M$ (or, equivalently, a gem $\Gamma \in G_s^{(4)}$ of $M$) and a cyclic permutation $\e=(\e_0,\e_1,\e_2,\e_3,4)$ of $\Delta_4$ exist, so that $\mathcal T(\Gamma,\e)=(H_0, H_1, H_2)$ is a gem-induced trisection of $M$, where $H_1$ and $H_2$ are 4-dimensional handlebodies and $H_0$ is a collar on $\partial M$; moreover, all pairwise intersections $H_0 \cap H_1, H_1\cap H_2, H_2 \cap H_0$ are $3$-dimensional handlebodies of  the same genus. 
In virtue of a celebrated theorem by \cite{Laudenbach-Poenaru}, together with its non-orientable version by \cite{Miller-Naylor}, a 4-dimensional handlebody of genus $m$ \, ${\overset{(\sim)}{\mathbb Y}}\!^4_m$ (orientable or  non-orientable, according to the orientability of $M$) may be glued in a unique way to the ``free'' boundary of $H_0$, so to obtain  $\bar H_0 = H_0 \cup {\overset{(\sim)}{\mathbb Y}}\!^4_m \cong {\overset{(\sim)}{\mathbb Y}}\!^4_m.$  It is now easy to check that the triple $(\bar H_0, H_1, H_2)$ actually constitutes a trisection of the closed 4-manifold $\bar M =  M \cup {\overset{(\sim)}{\mathbb Y}}\!^4_m$, with the same pairwise intersections as $\mathcal T(\Gamma,\e)$,  and hence with the same intersecting surface.     

The second part of the statement directly follows. 
\qed       

\medskip 

In the following, a trisection of a closed (orientable or non-orientable) 4-manifold $\bar M$ will be said to {\it arise from a colored triangulation}  if  it is either a gem-induced trisection of $\bar M$ or it is obtained from a gem-induced trisection of $M$ (such that $\bar M \cong M \cup {\overset{(\sim)}{\mathbb Y}}\!^4_m$) according to Theorem \ref{trisection_from_gem-induced}.     

Moreover, for sake of conciseness,  we will denote by $\mathbb S^{n-1}\otimes\mathbb S^1$ either $\mathbb S^{n-1}\times\mathbb S^1$ or $\mathbb S^{n-1}\tilde\times\mathbb S^1$, i.e. the orientable or non-orientable $\mathbb S^{n-1}$-bundle over $\mathbb S^1$.

\begin{example} \label{S3xS1}  {\em  
In \cite{Casali-Cristofori gem-induced}, 
$\mathbb Y^4_1$ has been proved to admit a gem-induced trisection of genus $1$, through its well-known order eight crystallization; on the other hand, Example \ref{example-handlebodies} proves that  $\tilde{\mathbb Y}^4_1$ admits a gem-induced trisection of genus $1$, too, through the (order eight) crystallization of Figure \ref{Y4_1-no}.  \    Hence,  \  \ 
$g_T(\mathbb S^3 \otimes \mathbb S^1)  \le 1$ follows from Theorem \ref{trisection_from_gem-induced}. Actually, the equality holds via the well-known characterization of $\mathbb S^4$ as the only closed $4$-manifold with trisection genus equal to zero (see also Proposition \ref{GTgenus-properties}(ii)).   
Thus,  both  $\mathbb S^3 \times \mathbb S^1$ and $\mathbb S^3 \tilde \times \mathbb S^1$ are proved to admit a  trisection with minimal genus arising from  a colored triangulation.  }
\end{example}

We want now to describe a possible way to pass from a gem of the closed 4-manifold  $\bar M$ to a gem of a compact $4$-manifold  $M$ so that $\bar M \cong M \cup {\overset{(\sim)}{\mathbb Y}}\!^4_m$; to this aim, the notion of $\rho$-pair\footnote{$\rho$-pairs and their switching were introduced by Lins (\cite{Lins-book}) and subsequently studied in \cite{Casali-Cristofori ElecJComb 2015}, \cite{Bandieri-Gagliardi}, \cite{CFMT}.} turns out to be very useful. 

\begin{definition}\label{rho-pair} {\em 
A $\rho_h${\it -pair} ($1\leq h\leq n$) of color $i\in\Delta_n$ in an $(n+1)$-colored graph $\Gamma$ is a pair of $i$-colored edges $(e,f)$ sharing the  same $\{i,c\}$-colored cycle for each $c\in\{c_1,\ldots,c_h\}\subseteq\Delta_n.$ Colors $c_1,\ldots,c_h$ are said to be {\it involved}, while the other $n-h$ colors are said to be {\it not involved} in the $\rho_h$-pair.
\\
The {\it switching} of $(e,f)$ consists in canceling $e$ and $f$ and establishing new $i$-colored edges between their endpoints; the reversed operation is obviously the switching of a $\rho_{n-h}$-pair.  
Although, in general, the switching may be performed in two different ways, it is uniquely determined if $\G \in G^{(n)}_s,\ h\in \{n-1, n\}$ and the bipartition of each non-singular $\hat c$-residue is preserved. } 
\end{definition} 

The topological effects of the switching of $\rho_{n-1}$- and $\rho_n$-pairs have been completely determined in the case of closed $n$-manifolds: see \cite{Bandieri-Gagliardi}, where it is proved that a $\rho_{n-1}$-pair (resp. $\rho_n$-pair) switching does not affect the represented $n$-manifold  (resp. either induce the splitting into two connected summands, or the ``loss'' of a $\mathbb S^{n-1} \otimes \mathbb S^1$  summand in the represented $n$-manifold). 

\begin{proposition}\label{trisection_from_rho1-pairs} 
Let $\bar M$ be a closed 4-manifold and $\Gamma$ a crystallization of $\bar M$. 
If  $\Gamma$ contains $m$ $\rho_1$-pairs of color $i$ ($i \in \Delta_3$) and involving color $4$, so that the sequence of $m$ switchings yields a 5-colored graph \ $\Gamma^\prime  \in G^{(4)}_s$ admitting a gem-induced trisection $\mathcal T (\Gamma^\prime, \varepsilon)$ (for a suitable cyclic permutation $\e=(\e_0,\e_1,\e_2,\e_3,4)$ of $\Delta_4$), 
then 
 $$g_T(\bar M) \le \rho_{\varepsilon_{\hat 4}}(\G_{\hat 4}) +m. $$ 
\end{proposition}
 
\begin{figure}[t]
\centering
\scalebox{0.42}{\includegraphics{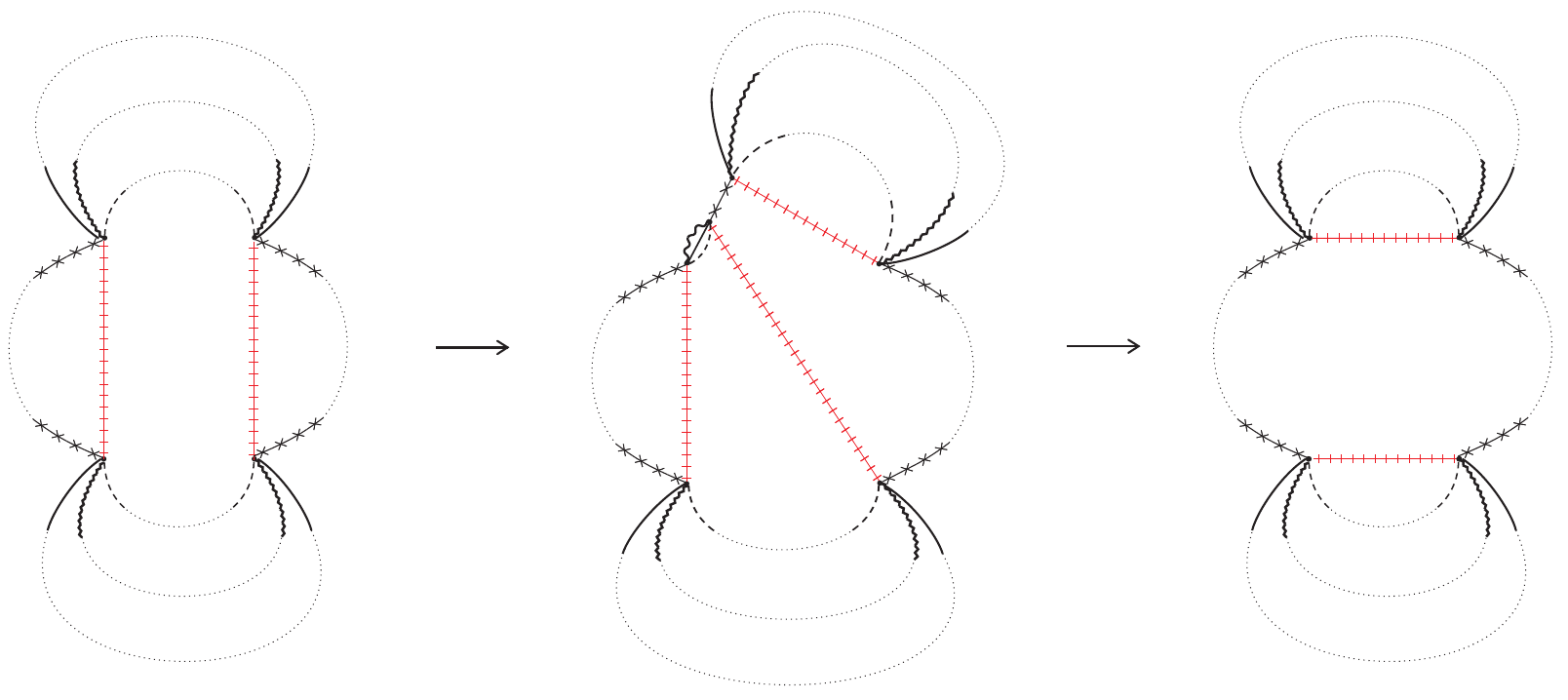}}
\caption{Factorization of a $\rho_3$-pair switching into an edge insertion and a 3-dipole cancellation}
\label{fig: ro3} 
\end{figure}

\dimo      
First, let us consider the case $m=1$: let $\Gamma^\prime \in G^{(4)}_s$ be the 5-colored graph obtained from $\Gamma$ (representing the closed $4$-manifold $\bar M$) by switching a $\rho_1$-pair of color $i$ ($i \in \Delta_3$) and involving color $4$. It is easy to check that $\Gamma^\prime$ represents a compact $4$-manifold  $M$ whose  
boundary   is  PL-homeomorphic to $\mathbb S^2 \otimes \mathbb S^1$: in fact, the $\hat 4$-residue $\Gamma^\prime_{\hat 4}$ represents  the orientable or non-orientable $\mathbb S^2$-bundle over $\mathbb S^1$ (according to  $\Gamma^\prime$ being bipartite or not), since it contains a $\rho_3$-pair whose switching yields the $\hat 4$-residue $\Gamma_{\hat 4}$ representing $\mathbb S^3$, while all other $\hat c$-residues $\Gamma^\prime_{\hat c}$ ($c \in \Delta_3$) represent $\mathbb S^3,$ since they give rise to $\Gamma_{\hat c},$ representing $\mathbb S^3$, by switching  a $\rho_2$-pair. 
Moreover, $\bar M$ is obtained from $M$ by attaching a $3$-handle to its boundary, and this $3$-handle is orientable or non-orientable according to the orientability of  $\partial M$.
In fact, as  shown  in Figure \ref{fig: ro3}, the switching of the $\rho_3$-pair in  $\Gamma^\prime$ can be factorized by inserting a $4$-colored edge (whose end-points belong to different bipartition classes in $\Gamma^\prime_{\hat 4}$ if and only if $\Gamma$ is bipartite) and subsequently cancelling a {\it $3$-dipole}, i.e. a subgraph  consisting of two vertices joined by three colored edges, so that the vertices belong to different bicolored cycles involving the remaining colors. 
It is not difficult to check that the insertion of the $4$-colored edge corresponds to  ``breaking''   a tetrahedral boundary face of $M$ and inserting a new pair of $4$-simplices sharing the same $3$-dimensional face opposite to the $4$-labelled vertex,  so to transform the boundary into a $3$-sphere; hence, this operation may be seen as the attachment of a polyhedron homeomorphic to $\mathbb D^3\times\mathbb D^1$  to the boundary of $M$, without affecting its interior\footnote{An analogous   argument has already been  used  in \cite[Proposition 11(ii)]{Casali-Cristofori Kirby-diagrams},  to prove the correspondence between the switching of $\rho_3$-pairs in a graph representing $[\#_m(\mathbb S^2 \times \mathbb S^1)] \times I$ and the attachment of $3$-handles.}. On the other hand,  the elimination of the $3$-dipole does not affect the represented manifold, since it corresponds to a re-triangulation of a subcomplex of $\vert K(\G)\vert$  homeomorphic to a $4$-ball (see \cite{generalized-genus} and references within). 

The case $m  > 1$ may be proved by induction via similar arguments, by taking into account that, at the $k$-th step ($2\leq k\leq m$), the switching of the considered $\rho_3$-pair transforms the boundary, represented by the $\hat 4$-residue, from $\#_k (\mathbb S^2\otimes\mathbb S^1)$ into $\#_{k-1}(\mathbb S^2\otimes\mathbb S^1).$   
Hence, $ \Gamma^\prime$ represents a compact 4-manifold $M$ such that $\bar M \cong M \cup {\overset{(\sim)}{\mathbb Y}}\!^4_m$.     

\smallskip

Now, if $\Gamma^\prime \in G^{(4)}_s$ admits a gem-induced trisection $\mathcal T (\Gamma^\prime, \varepsilon)$ (for a suitable cyclic permutation $\e=(\e_0,\e_1,\e_2,\e_3,4)$ of $\Delta_4$), $g_{GT} (M) \le  \rho_{\varepsilon_{\hat 4}}(\G^\prime_{\hat 4})$ holds by  Proposition \ref{GTgenus-properties}(i).    
On the other hand, it is easy to check that, by switching   a single $\rho_1$-pair of color $i$ involving color $4$, the number of  $\{i,j\}$-cycles, with $j \in \Delta_3,$ is decreased by one, while the number of  $\{r,s\}$-cycles, with $r,s \in \Delta_3-\{i\}$, is unaffected; as a consequence, by formula \eqref{eq.regulargenus} in Section 2, $\rho_{\varepsilon_{\hat 4}}(\G^\prime_{\hat 4})= \rho_{\varepsilon_{\hat 4}}(\G_{\hat 4}) + m$ holds.  

Hence, the thesis follows from Theorem \ref{trisection_from_gem-induced}:  
$$g_{T} (\bar M) \le g_{GT} (M) \le \rho_{\varepsilon_{\hat 4}}(\G_{\hat 4}) + m.  \eqno
$$ 
\qed      

\begin{figure}[t]
\centering
\scalebox{0.45}{\includegraphics{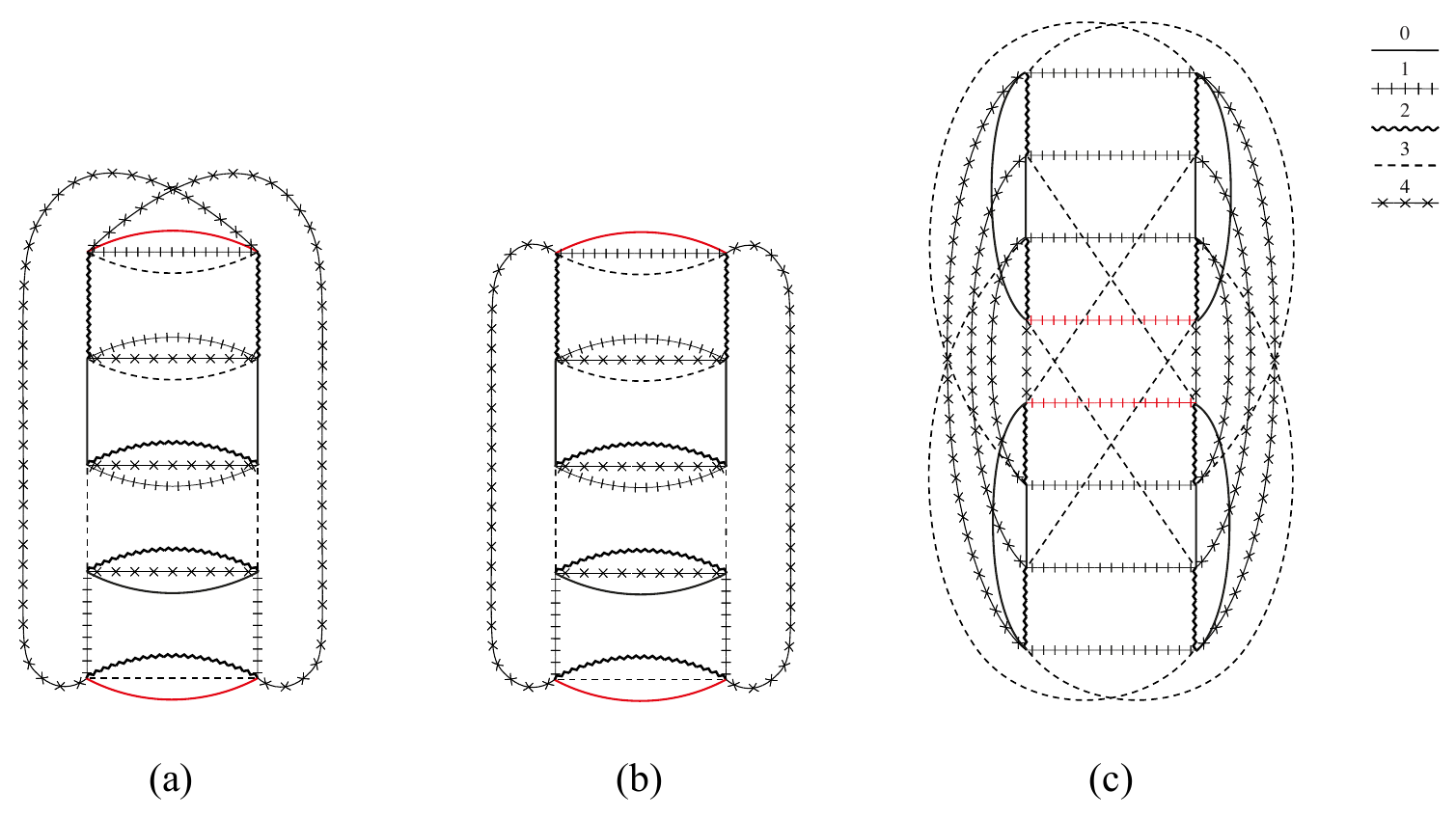}}
\caption{minimal crystallizations of  $\mathbb S^3 \times \mathbb S^1$,  $\mathbb S^3 \tilde \times \mathbb S^1$ and $\mathbb R \mathbb P^4$, with $\rho_1$-pairs in red}   
\label{fig.examples_trisections_from_gem-induced}
\end{figure} 

 \begin{example} \label{S3xS1_bis} {\em 
 $ g_T(\mathbb S^3 \otimes \mathbb S^1) =1$ may be re-obtained also via Proposition \ref{trisection_from_rho1-pairs}, starting from the crystallizations of $\mathbb S^3\times\mathbb S^1$ and $\mathbb S^3 \tilde \times \mathbb S^1$ depicted in Figure \ref{fig.examples_trisections_from_gem-induced}(a) and 4(b) respectively: in fact, it is easy to check that each of them admits  a $\rho_1$-pair of color $0$ involving color $4$, whose switching yields  a gem of $  {\overset{(\sim)}{\mathbb Y}}\!^4_1$ satisfying the condition of Proposition  \ref{CS gem-induced trisections}. Since both   starting crystallizations have $\rho_{\varepsilon_{\hat i}}(\G_{\hat i})=0$ for each $i \in \Delta_4$,   $ g_T(\mathbb S^3 \otimes \mathbb S^1) \le 0+1=1$ follows from Proposition \ref{trisection_from_rho1-pairs}. 
As already pointed out, the equality actually holds, since trisection genus zero characterizes $\mathbb S^4$ among all closed $4$-manifolds.  }
 \end{example}  
 
\begin{example} \label{RP4}   {\em 
$ g_T(\mathbb R \mathbb P^4) =  2$  is realized by a trisection arising from  a colored triangulation, too. In fact, it is easy to check that the crystallization of $\mathbb R \mathbb P^4$ depicted in Figure \ref{fig.examples_trisections_from_gem-induced}(c)  admits  a $\rho_1$-pair of color $1$ involving color $4$, whose switching yields a gem satisfying the condition of Proposition  \ref{CS gem-induced trisections}. Since the starting crystallization
has $\rho_{\varepsilon_{\hat i}}(\G_{\hat i})=1$ for each $i \in \Delta_4$,   $ g_T(\mathbb R \mathbb P^4) \le  1+1=2$ follows from Proposition \ref{trisection_from_rho1-pairs}; then the above claim is a direct consequence of inequality (\ref{genus_inequality}). } 
\end{example}   

\begin{example} \label{RP2xS2}  {\em 
$ g_T(\mathbb S^2 \times \mathbb R \mathbb P^2) =  3$ is realized by a trisection  arising from  a colored triangulation, too. 
In fact, it is easy to check that the  (order 24) crystallization of $\mathbb S^2 \times \mathbb R \mathbb P^2$ depicted in \cite[Figure 3]{Basak-Casali} admits a $\rho_1$-pair of color $1$ and involving color $0$; by a suitable permutation of colors, the $\rho_1$-pair of color $1$ involves color $4$, and its switching yields a gem satisfying the condition of Proposition  \ref{CS gem-induced trisections}. Since the above crystallization of $\mathbb S^2 \times \mathbb R \mathbb P^2$ has $\rho_{\varepsilon_{\hat i}}(\G_{\hat i})=2$ for each $i \in \Delta_4$,   $ g_T(\mathbb S^2 \times \mathbb R \mathbb P^2) \le  2+1=3$ follows from Proposition \ref{trisection_from_rho1-pairs}; then the  above claim is a direct consequence of inequality (\ref{genus_inequality}). } 
\end{example}

As already mentioned in the Introduction, the existence of  trisections of minimal genus arising from colored triangulations can now be proved for a large class of closed $4$-manifolds:

\begin{proposition} \label{calculations} Let \ 
$\bar M\, \cong_{PL}\,(\#_p\mathbb {CP}^2)\,\#\,(\#_{p^{\prime}}(-\mathbb {CP}^2))\, \#\, (\#_q(\mathbb S^2\times \mathbb S^2))\,\#$ $(\#_r(\mathbb S^3 \otimes \mathbb S^1))\, \#\, (\#_s\mathbb R \mathbb P^4)\,\# \, (\#_t K3),$  with $p,p^{\prime},q, r,s,t \geq 0$.  
Then, its trisection genus \  $ g_T(\bar M)  = (p + p^{\prime} + 2q + 22t)+r+2s$ \ 
is realized by a trisection arising from a colored triangulation.
\end{proposition} 

\dimo       
For each summand, the existence of a trisection of minimal genus arising from a colored  triangulation (or, equivalently, from a suitable gem) is ensured by Examples \ref{S3xS1} (or  \ref{S3xS1_bis}) and \ref{RP4}, together with results in \cite{SpreerTillmann(Exp)}  concerning minimal gem-induced trisections of $\mathbb C \mathbb P^2$, $\mathbb S^2 \times \mathbb S^2$ and the K3-surface.  
Note that all the involved gems represent closed $4$-manifolds if and only if  $r+s =0$. 

Now, the so called {\it graph connected sum} may be performed on the above gems, yielding a gem of the compact 4-manifold, with empty or connected boundary, obtained by connected sum (or by boundary connected sum, if  $r+s\ge 2$) of the represented manifolds. 
Indeed, according to gem theory, given two ($n+1$)-colored graphs $\G_1$ and $\G_2$,  their graph connected sum with respect to vertices $v_1$ and $v_2$  ($v_i$ in $\G_i$, $\forall i\in \{1,2\}$) is the graph $\G_1 \#_{v_1,v_2} \G_2$  obtained by deleting $v_1$ and $v_2$ and by welding the hanging edges of the same color;  if $\G_1$ (resp. $\G_2$) represents the $n$-manifold $M_1$ (resp. $M_2$), then  $\G_1 \#_{v_1,v_2} \G_2$ is known to represent $M_1 \# M_2$  in case at least one of $M_1$ and $M_2$ is closed, and  $M_1 \natural M_2$  in case both $M_1$ and $M_2$ have connected non-empty  boundary  (see  \cite[Section 7]{Grasselli-Mulazzani} for details).   

Moreover, it is not difficult to check that, by suitably performing the above sequence of graph connected sums, the obtained $5$-colored graph yields a gem-induced trisection with genus equal to the sum of the genera of the  trisections (assumed to be of minimal genus) of each summand: see Proposition \ref{GTgenus-properties}(iv) and the proof of \cite[Proposition 26(ii)]{Casali-Cristofori gem-induced}.  

The existence of the required trisection of $\bar M$ follows by making use of  Theorem \ref{trisection_from_gem-induced},  in case $r+s >0$.   

\smallskip	
On the other hand, the trisection genus of $\bar M$ exactly coincides with the sum of the genera, in virtue of inequality   \eqref{genus_inequality}.  
\qed

\medskip

Let us now restrict the attention to the orientable case, where trisections arising from colored triangulations can be proved to exist for each closed $4$-manifold, via Kirby diagrams and associated colored graphs.  

In fact, if  $\bar M$ is a closed orientable $4$-manifold, it is well-known that - in virtue of the already cited  result by Laudenbach-Poenaru (\cite{Laudenbach-Poenaru}) -   a Kirby diagram $(L,d)$ of $\bar M$ describes only the attachments of the 1-handles (via 
dotted components) and 2-handles (via framed components) of a handle decomposition of $\bar M$. Hence, $(L,d)$ actually coincides with a Kirby diagram of the compact 4-manifold $M$ consisting only of the $h$-handles of $\bar M$, with $h \le 2$. Obviously, if no $3$-handle appears, $M$  can be identified with $\bar M$, up to capping off its spherical boundary. 

On the other hand, by results in \cite{Casali-Cristofori Kirby-diagrams} and \cite[Section 4]{Casali-Cristofori gem-induced}, a gem-induced trisection of $M$ can be algorithmically constructed starting from a Kirby diagram.  Thus, Theorem \ref{trisection_from_gem-induced} directly yields the required trisection of $\bar M$ arising from a colored triangulation. 

Let us specify that, as in \cite{Casali-Cristofori Kirby-diagrams} and \cite{Casali-Cristofori gem-induced}, the considered Kirby diagrams are {\it connected}, in the sense that their associated planar projections are supposed to be connected; moreover, %.Moreover, 
the (possible) dotted components of $(L,d)$ are assumed to be in {\it good position}, i.e. they are unknotted, unlinked and with overcrossings and undercrossings never alternating along them. Note that - without loss of generality - any closed orientable 4-manifold can be represented by a Kirby diagram satisfying  these conditions (see \cite{GS}).

\bigskip 

According to Theorem \ref{trisection_from_Kirby-diagram} stated in the Introduction, we are now able to extend to the whole class of closed orientable 4-manifolds the estimation of the trisection genus via Kirby diagrams that was already obtained in  \cite{Casali-Cristofori gem-induced} under the assumption of the existence of a handle decomposition lacking in $3$-handles.

\bigskip 

{\it Proof of Theorem \ref{trisection_from_Kirby-diagram}. }
Let $\bar M$ be a closed orientable $4$-manifold represented by a (connected) Kirby diagram $(L,d)$ whose dotted components, if any, are in good position.  

Further, let us suppose that the handle decomposition of $\bar M$ associated to $(L,d)$ contains  $q>0$ 3-handles (if $q=0,$ both estimations were already proved in \cite[Corollary 4]{Casali-Cristofori gem-induced}) and let $M$ be the compact $4$-manifold consisting only of the $h$-handles, with $h \le 2$, of the above handle decomposition. 

In the general case when $(L,d)$ contains both dotted and framed components (i.e. the handle decomposition of $M$ contains both $1$-handles and $2$-handles: case (i) of the statement), \cite[Theorem 12]{Casali-Cristofori Kirby-diagrams}  and \cite[Theorem 3]{Casali-Cristofori gem-induced} allow to construct a  gem of $M$, so that it admits a genus $s+1$  gem-induced trisection ($s$ being the crossing number of $(L,d)$).  
On the other hand, if  $(L,d)$ contains no dotted components (i.e. the handle decomposition of $M$ contains only $2$-handles: case (ii) of the statement), \cite[Theorem 7(ii)]{Casali-Cristofori Kirby-diagrams}  and \cite[Theorem 3]{Casali-Cristofori gem-induced} yield  again a gem of $M,$ so that it admits  
a genus $m_{\alpha}$ gem-induced trisection ($m_{\alpha}$ being the number of $\alpha$-colored regions in a chess-board coloration of $(L,d)$, where $\alpha$ is the color of the unbounded region). 

In both cases, the thesis follows from Theorem \ref{trisection_from_gem-induced}, since $\bar M \cong M \cup \mathbb Y^4_q$: in case (i) (resp. (ii)), we have
$$ g_T(\bar M) \le g_{GT} (M) \le s+1  \ \  \  (\text{resp. }  g_T(\bar M) \le g_{GT} (M) \le  m_\alpha).$$

\vskip-0.3truecm 
\ \qed 
\medskip

\begin{example} \label{S3xS1_ter}  {\em 
$ g_T(\mathbb S^3 \times \mathbb S^1) =1$ may be re-obtained also via Theorem \ref{trisection_from_Kirby-diagram}(i), starting from the Kirby diagram of $\mathbb S^3 \times \mathbb S^1$ consisting only of one dotted component: in fact, since there are no crossings, $ g_T(\mathbb S^3 \times \mathbb S^1) \le 0+1=1$ trivially follows.  }
\end{example}

\begin{remark}\label{dual-handle-decomposition} {\em 
We point out that an estimation of the trisection genus similar to the one given in Theorem \ref{trisection_from_Kirby-diagram}(i) has been obtained in \cite{Kepplinger} via trisection diagrams.

\noindent
Note also that, when there are no dotted components, the inequality in Theorem \ref{trisection_from_Kirby-diagram}(ii) should be preferred to the one in  Theorem \ref{trisection_from_Kirby-diagram}(i), since it obviously gives a better estimation.  
Therefore, for example, if $\bar M$ is a closed orientable $4$-manifold admitting a handle decomposition with no $3$-handles, it  could be more convenient to apply statement (ii) to the Kirby diagram representing the dual handle decomposition. }   
\end{remark}

\bigskip

We conclude the paper by pointing out that Theorem \ref{trisection_from_Kirby-diagram} can be easily generalized to the case of disconnected Kirby diagrams. In fact, it is well-known that, in this case, the represented compact manifold is the boundary connected sum of the compact manifolds represented by each connected component; hence, the estimation of the trisection genus directly follows from multiple application of Theorem \ref{trisection_from_Kirby-diagram}, together with subadditivity (see Proposition \ref{GTgenus-properties}(iv) and the proof of  Proposition \ref{calculations}):

\begin{corollary}\label{disconnected_diagrams}
Let  $\bar M$ be a closed orientable $4$-manifold and $(L,d)$ a Kirby diagram of $\bar M$ with $c$ connected components and whose dotted components - if any - are in good position. Then:
\begin{itemize}
\item[(i)]
 $$g_T(\bar M) \le s+c ,$$  $s$ being the crossing number of $(L,d)$. 
\item[(ii)]
Furthermore, if  $(L,d)$ has no dotted components,
then 
$$g_T(\bar M) \le m_\alpha+c-1,$$
$m_{\alpha}$ being  the number of $\alpha$-colored regions in a chess-board coloration of $(L,d)$.
\end{itemize}
\end{corollary}
\vskip-0.3truecm 
\ \qed 
\bigskip

\noindent {\bf Acknowledgements:\ }  This work was supported by GNSAGA of INDAM and by the University of Modena and Reggio Emilia, project:  {\it ``Discrete Methods in Combinatorial Geometry and Geometric Topology''}.

\end{document}